%% file: main.tex
\documentclass[leqno]{amsart}

\usepackage[foot]{amsaddr}

\usepackage{color}
\usepackage{geometry}
\usepackage{subfiles}
\usepackage{graphicx}
\usepackage[section]{placeins}
\usepackage{hyperref}
\usepackage{stmaryrd}
\usepackage{amsmath,amssymb,amsfonts,amsthm,textcomp}
\usepackage{mathtools}
\usepackage{tensor}
\usepackage{multicol}
\usepackage{subfig}
\usepackage{tikz}
\usetikzlibrary{arrows}
\usetikzlibrary{calc}
\graphicspath{
              {./figures_section-1/}
              {./figures_section-2/}
              {./figures_section-3/}
             }

\input{commands}

\begin{document}

\title[On the equivalence of different adaptive batch size selection strategies for SGD]{On the equivalence of different adaptive batch size selection strategies for stochastic gradient descent methods}
\author{Luis Espath$^{1,\sharp}$, Sebastian Krumscheid$^2$, Ra\'{u}l Tempone$^{3,4,5}$ \& Pedro Vilanova$^6$}
\address{$^1$School of Mathematical Sciences, University of Nottingham, Nottingham, NG7 2RD, United Kingdom}
\address{$^2$Department of Mathematics, Karlsruhe Institute of Technology, 76131 Karlsruhe, Germany.}
\address{$^3$Department of Mathematics, RWTH Aachen University, Geb\"{a}ude-1953 1.OG, Pontdriesch 14-16, 161, 52062 Aachen, Germany.}
\address{$^4$King Abdullah University of Science \& Technology (KAUST), Computer, Electrical and Mathematical Sciences \& Engineering Division (CEMSE), Thuwal 23955-6900, Saudi Arabia.}
\address{$^5$Alexander von Humboldt Professor in Mathematics for Uncertainty Quantification, RWTH Aachen University, Germany.}
\address{$^6$Department of Mathematical Sciences, Stevens Institute of Technology, Hoboken, NJ 07030 USA.}
\email{$^\sharp$ espath@gmail.com}

\date{\today}

\begin{abstract}
\noindent
\subfile{./abstract.tex}
\end{abstract}

\maketitle

\tableofcontents                        


\subfile{./section-1.tex}

\subfile{./section-2.tex}

\subfile{./section-3.tex}

\subfile{./section-4.tex}


%


\footnotesize


\end{document}

%% file: commands.tex
\newfont{\tenbbb}{msbm10}
\newfont{\svnbbb}{msbm8}

\newcommand{\bs}[1]{\boldsymbol{#1}}
\newcommand{\cl}[1]{\mathcal{#1}}

\newcommand{\bb}[1]{\mathbb{#1}}

\newcommand{\id}{{1\!\!1}}

\newcommand{\norm}[2]{\left\lVert {#1} \right\rVert_{#2}}

\newcommand{\tr}{\mathrm{tr}\mskip2mu}

\newtheorem{lem}{Lemma}
\newtheorem{prop}{Proposition}

\newtheorem{set}{Definition}
\newtheorem{ass}{Assumption}


%% file: abstract.tex
In this study, we demonstrate that the norm test and inner product/orthogonality test presented in \cite{Bol18} are equivalent in terms of the convergence rates associated with Stochastic Gradient Descent (SGD) methods if $\epsilon^2=\theta^2+\nu^2$ with specific choices of $\theta$ and $\nu$. Here, $\epsilon$ controls the relative statistical error of the norm of the gradient while $\theta$ and $\nu$ control the relative statistical error of the gradient in the direction of the gradient and in the direction orthogonal to the gradient, respectively. Furthermore, we demonstrate that the inner product/orthogonality test can be as inexpensive as the norm test in the best case scenario if $\theta$ and $\nu$ are optimally selected, but the inner product/orthogonality test will never be more computationally affordable than the norm test. Finally, we present two stochastic optimization problems to illustrate our results.
\\
\noindent\textbf{AMS subject classifications:}
$\cdot$
65K05 
$\cdot$
90C15 
$\cdot$
65C05 
$\cdot$

%% file: section-1.tex
\section{Introduction}

The selection of a suitable sample size for estimating the expected gradient is a common aspect of approximations to stochastic optimization problems
using Stochastic Gradient Descent (\texttt{SGD}) methods. Indeed, the
sample size is
linked to
the statistical estimator's desired accuracy.
In \cite{Car91,Byr12,Has14,Car18}, the statistical error is defined in
terms of the expected value of the squared norm of the difference
between the estimated mean gradient
and the exact expected gradient. Using this definition of the statistical error, the norm test provides a criterion for determining the sample size. Conversely, in a seminal work, Bollapragada et al.\
\cite{Bol18}, among many contributions, propose the inner
product/orthogonality test (inner/orth for short) in lieu of the norm test to determine the sample size. The inner/orth test decomposes the statistical error in two directions. The inner product test is performed in reference to the parallel direction, that is, the
direction of the exact mean gradient, whereas the orthogonal test is performed in relation to the normal component.

To state the stochastic optimization problem, let
$\bs{\xi}\in \bb{R}^{d_{\bs{\xi}}}$ be the design variable in
dimension $d_{\bs{\xi}}\in\bb{N}$ and
$\bs{\vartheta}\in\bb{R}^{d_{\bs{\theta}}}$ be a vector-valued random
variable in dimension $d_{\bs{\theta}}\in\bb{N}$, whose probability
distribution $\pi$ may depend on $\bs{\xi}$. We
assume throughout this work that we can generate as many independent and identically
distributed (i.i.d.) samples from $\pi$ as we require. Here,
$\bb{E}[\cdot|\bs{\xi}]$, $\bb{V}[\cdot|\bs{\xi}]$, and
$\bb{C}[\cdot,\cdot|\bs{\xi}]$ are the expectation,
variance, and covariance operators conditioned on $\bs{\xi}$, respectively. Aiming
to optimize conditional expectations on $\bs{\xi}$, we state our
problem as follows,
\begin{equation}\label{eq:underlying.problem}
\bs{\xi}^* = \underset{\bs{\xi} \in \bb{R}^{d_{\bs{\xi}}}}{\arg\min} \, \bb{E} [f(\bs{\xi},\bs{\vartheta})|\bs{\xi}],
\end{equation}
where
$f \colon \bb{R}^{d_{\bs{\xi}}} \times \bb{R}^{d_{\bs{\theta}}} \to \bb{R}$
is given. Let the objective function in our
problem be denoted by
$F(\bs{\xi})\coloneqq\bb{E} [f(\bs{\xi},\bs{\vartheta})|\bs{\xi}]$. In
minimizing \eqref{eq:underlying.problem} with respect to the design
variable $\bs{\xi} \in \bb{R}^{d_{\bs{\xi}}}$, \texttt{SGD} is
obtained using the following updating rule
\begin{equation}\label{eq:sgd}
\bs{\xi}_{k+1} = \bs{\xi}_{k} - \eta_k \bs{\upsilon}_k\;,\quad k\in\bb{N}_0\;,
\end{equation}
where $\bs{\upsilon}_k$ is an estimator of the expected gradient at
$\bs{\xi}_{k}$, i.e., $\bs{\upsilon}_k\approx \nabla_{\bs{\xi}} F(\bs{\xi}_{k})$
in a statistical sense made precise below.

In the norm test, the number of samples required to control the relative
statistical error of the mean gradient estimator $\bs{\upsilon}_k$ is
identified in terms of a prescribed a maximum tolerance, namely
$\epsilon>0$. Conversely, in the inner/orth test, the relative
statistical error of the estimated mean gradient is decomposed into an
error in the direction of the exact gradient
and its orthogonal component.
Each of these errors is respectively controlled by a tolerance
$\theta\ge0$ and $\nu\ge0$. In this study, we demonstrate that if the
inner/orth test criterion is met, then the norm test will be
automatically satisfied for all
$\epsilon^2\ge\theta^2+\nu^2$. Furthermore, these two approaches are
equivalent in terms of computational cost and convergence rates as
long as
\begin{equation*}
\epsilon^2=\theta^2+\nu^2,
\end{equation*}
with $\theta$ and $\nu$ satisfying the condition
\begin{equation*}
\dfrac{\theta^2}{\nu^2}=\dfrac{\bs{\Sigma}\colon\bs{P}_{\!\scriptscriptstyle{\nabla}}}{\bs{\Sigma}\colon\bs{P}_{\!\scriptscriptstyle{\perp}}}\;.
\end{equation*}
Here, $\bs{\Sigma}$ denotes the covariance of the gradient estimator,
$\bs{P}_{\!\scriptscriptstyle{\nabla}}\coloneqq\bs{e}_{\scriptscriptstyle{\nabla}}\otimes\bs{e}_{\scriptscriptstyle{\nabla}}$,
$\bs{P}_{\!\scriptscriptstyle{\perp}}\coloneqq\id-\bs{e}_{\scriptscriptstyle{\nabla}}\otimes\bs{e}_{\scriptscriptstyle{\nabla}}$,
and
$\bs{e}_{\scriptscriptstyle{\nabla}}\coloneqq\nabla_{\bs{\xi}}F(\bs{\xi})/\norm{\nabla_{\bs{\xi}}
  F(\bs{\xi})}{}$, where we have suppressed the dependence on the
iteration $k$ for brevity.

The rest of this study is structured as follows. In \S\ref{fundations}, we present assumptions and definitions. In \S\ref{stat.decomposition}, we present a Lemma demonstrating the equivalence of the inner product/orthogonality test and the norm test for determining the sample size. In \S\ref{cov.decomposition}, we demonstrate how the stochastic gradient's covariance is decomposed to compute the inner product/orthogonality test to determine the sample size. In \S\ref{complexity}, we compare the complexity equivalence of the inner product/orthogonality test and the norm test in terms of the sample size. In \S\ref{numerics}, we show the numerical predictions for both tests.

\section{Foundations}
\label{fundations}

Throughout this study, we use $\norm{\cdot}{}$ to denote the usual
Euclidean norm. Furthermore, we consider the following assumptions to
hold for the stochastic optimization
problem~\eqref{eq:underlying.problem} and its \texttt{SGD}
approximation~\eqref{eq:sgd}.
\begin{ass}[$L$-Lipschitz gradient ($L$-smoothness)]\label{as:lipschitz}
  The gradient of $F \colon \bb{R}^{d_{\bs{\xi}}}\mapsto\bb{R}$ is
  $L$-Lipschitz on the feasible set
  $\Xi\subseteq \bb{R}^{d_{\bs{\xi}}}$ for some $L > 0$, in the sense
  that
\begin{equation}\label{eq:lipschitz}
\norm{\nabla_{\bs{\xi}} F(\bs{x}) - \nabla_{\bs{\xi}} F(\bs{y})}{} \le L \norm{\bs{y} - \bs{x}}{}, \qquad \forall \bs{x}, \bs{y} \in \Xi\;.
\end{equation}
\end{ass}
\begin{ass}[Convex and Strongly convex]\label{as:strongly.convex}
  The objective function $F$ is $L$-smooth convex and $\mu$-strongly
  convex such that for some $\mu>0$ it holds that
  \begin{equation}\label{eq:convex}
    F(\bs{y}) \le F(\bs{x}) + \langle \nabla_{\bs{\xi}} F(\bs{x}), \bs{y} - \bs{x} \rangle + \dfrac{L}{2} \norm{\bs{y} - \bs{x}}{}^2, \qquad \forall\,\bs{x}, \bs{y} \in \Xi\;,
  \end{equation}
  as well as
  \begin{equation}\label{eq:strongly.convex}
    F(\bs{y}) \ge F(\bs{x}) + \langle \nabla_{\bs{\xi}} F(\bs{x}), \bs{y} - \bs{x} \rangle + \dfrac{\mu}{2} \norm{\bs{y} - \bs{x}}{}^2, \qquad \forall\,\bs{x}, \bs{y} \in \Xi\;.
  \end{equation}
\end{ass}
\begin{ass}[Unbiased estimator]\label{as:grad_ass_unbiased}
  The estimated mean gradient $\bs{\upsilon}_k$ is a conditionally
  unbiased estimator of the exact gradient $\nabla_{\bs{\xi}} F$ at $\bs{\xi}_k$ for all $k\ge 0$ such that
\begin{equation}\label{eq:grad_ass_unbiased}
\bb{E}[\bs{\upsilon}_k |\bs{\xi}_k] = \nabla_{\bs{\xi}} F(\bs{\xi}_k)\;,\qquad\forall\,k\in\bb{N}_0\;.
\end{equation}
\end{ass}

Next, we introduce the norm test criterion for~\eqref{eq:sgd}.
\begin{set}[Norm test]\label{df:grad_rel_err.norm}
  Let $\bs{\upsilon}_k$ be a conditionally unbiased mean gradient
  estimator satisfying Assumption~\ref{as:grad_ass_unbiased} for all
  $k\in\bb{N}_0$. The statistical estimator $\bs{\upsilon}_k$ is said
  to satisfy the norm test for a given tolerance $\epsilon>0$, iff
\begin{equation}\label{eq:grad_rel_err.norm}
\bb{E}[\norm{\bs{\upsilon}_k-\nabla_{\bs{\xi}} F(\bs{\xi}_k)}{}^2|\bs{\xi}_k]
\le\norm{\nabla_{\bs{\xi}} F(\bs{\xi}_k)}{}^2\epsilon^2\;,\quad k\in\bb{N}_0\;.
\end{equation}
\end{set}

Using the norm test, we have the following characterization of the
\texttt{SGD}'s accuracy follows from the results in~\cite{Car21}.
\begin{prop}[Proposition 1 in \cite{Car21}]\label{prp:convergence.mice}
  Assume the mean gradient estimators $\bs{\upsilon}_k$ satisfy the
  norm test \eqref{eq:grad_rel_err.norm} for certain relative statistical
  tolerance $\epsilon>0$. If Assumptions \ref{as:lipschitz},
  \ref{as:strongly.convex}, and \ref{as:grad_ass_unbiased} hold and
  constant uniform step-size
  \begin{equation}\label{eq:step.size.strongly.convex}
    \eta_k\equiv\eta = \dfrac{2}{(L+\mu)(1 + \epsilon^2)}
  \end{equation}
  is used in the \texttt{SGD} method~\eqref{eq:sgd}, then
  the \texttt{SGD} method enjoys a linear convergence rate, namely
\begin{align}\label{eq:residual.recursion}
\bb{E}[\norm{\bs{\xi}_{k+1} - \bs{\xi}^*}{}^2] &\le \left[\frac{\left(\dfrac{\kappa-1}{\kappa+1}\right)^2+\epsilon^2}{1+\epsilon^2}\right]^k\bb{E}[\norm{\bs{\xi}_0 - \bs{\xi}^*}{}^2],
\end{align}
where $L$ is the Lipschitz gradient constant $L$ given in
\eqref{eq:lipschitz}, $\mu$ the strongly convex constant in
\eqref{eq:strongly.convex}, and $\kappa=L/\mu$ the condition number.
\end{prop}
Finally, note that, using the Lipschitz condition and
convexity, expression \eqref{eq:residual.recursion} may be written in
terms of the optimality gap, that is, in terms of
$F(\bs{\xi}_k)-F(\bs{\xi}^\ast)$.

\subsection{Statistical error decomposition}
\label{stat.decomposition}

In this section, we decompose the statistical error of the mean
gradient estimator into two components. These components are the
statistical error in the orthogonal and parallel directions
of the true gradient $\nabla_{\bs{\xi}}F$. Based on this decomposition, we
demonstrate the equivalence of two
different strategies commonly used to estimate the sample size of the
gradient estimator in terms of convergence rates.
For clarity, we will suppress the dependence of the statistical
estimators $\bs{\upsilon}_k$ and the evaluation points $\bs{\xi}_k$ in the following.

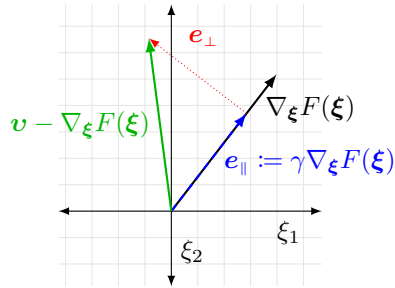
\begin{figure}[!h]
\centering
\begin{tikzpicture}[line/.style={>=latex}]
\coordinate (V1) at (1, 1.3);
\coordinate (V2) at (-1.3, 1.0);
\coordinate (V3) at ($(V1) + (V2)$);
\coordinate (V4) at ($1.4*(V1)$);
\draw[step=10pt, color=black!10] (-1.5, -1) grid (2, 2.75);
\draw[<->, line] (-1.5, 0) -- node [below, very near end] {$\xi_1$} (2, 0);
\draw[<->, line] (0, -1) -- node [right, very near start] {$\xi_2$} (0, 2.75);
\draw[->, line, color=black, thick] (0, 0) -- node [right=2pt, near end]    {$\nabla_{\bs{\xi}}F(\bs{\xi})$} (V4);
\draw[->, line, color=blue, thick, densely dashed] (0, 0) -- node [right=2pt] {$\bs{e}_{\scriptscriptstyle{\parallel}}\coloneqq\gamma\nabla_{\bs{\xi}}F(\bs{\xi})$} (V1);
\draw[->, line, color=red, densely dotted] (V1) -- +(V2) node [left=-30pt] {$\bs{e}_{\scriptscriptstyle{\perp}}$};
\draw[->, line, color=green!70!black, thick] (0, 0) -- node [left] {$\bs{\upsilon}-\nabla_{\bs{\xi}}F(\bs{\xi})$} (V3);
\end{tikzpicture}
\caption{Statistical error decomposition.}
\label{fg:error}
\end{figure}

The error decomposition presented in
Lemma~\ref{lm:error.decomposition} below is based on the idea illustrated
in Figure~\ref{fg:error}, which presents a parallel/orthogonal
decomposition of the error associated with the usage of a statistical estimator
$\bs{\upsilon}$. In this figure,
$\bs{\upsilon}-\nabla_{\bs{\xi}}F(\bs{\xi})=\bs{e}_{\scriptscriptstyle{\perp}}+\bs{e}_{\scriptscriptstyle{\parallel}}$
where
$\bs{e}_{\scriptscriptstyle{\perp}}\coloneqq\bs{\upsilon}-(1+\gamma)\nabla_{\bs{\xi}}F(\bs{\xi})$
and
$\bs{e}_{\scriptscriptstyle{\parallel}}\coloneqq\gamma\nabla_{\bs{\xi}}F(\bs{\xi})$
for an admissible $\gamma$ such that
$\bs{e}_{\scriptscriptstyle{\perp}}\cdot\bs{e}_{\scriptscriptstyle{\parallel}}=0$.
Indeed, the proof of the following result explicitly identifies $\gamma$, in addition to
$\bs{e}_{\scriptscriptstyle{\perp}}$ and $\bs{e}_{\scriptscriptstyle{\parallel}}$.
\begin{lem}{(Statistical error decomposition)}\label{lm:error.decomposition}
  Suppose Assumption~\ref{as:grad_ass_unbiased} holds for $\bs{\upsilon}$.
  Then the following statistical error decomposition holds
  \begin{equation}\label{eq:stat.decomposition}
    \bb{E}[\norm{\bs{\upsilon}-\nabla_{\bs{\xi}}F(\bs{\xi})}{}^2|\bs{\xi}]=\bb{V}[\bs{\upsilon}\cdot\bs{e}_{\scriptscriptstyle{\nabla}}|\bs{\xi}]+\bb{E}[\norm{\bs{\upsilon} - (\bs{\upsilon}\cdot\bs{e}_{\scriptscriptstyle{\nabla}})\bs{e}_{\scriptscriptstyle{\nabla}}}{}^2|\bs{\xi}],
  \end{equation}
  where $\bs{e}_{\scriptscriptstyle{\nabla}}\coloneqq\nabla_{\bs{\xi}}F(\bs{\xi})/\norm{\nabla_{\bs{\xi}} F(\bs{\xi})}{}$.
\end{lem}

\begin{proof}
By orthogonality, we have that
\begin{align}
0={}&\bs{e}_{\scriptscriptstyle{\perp}}\cdot\nabla_{\bs{\xi}}F(\bs{\xi}),\nonumber\\[4pt]
={}&\bs{\upsilon}\cdot\nabla_{\bs{\xi}}F(\bs{\xi})-(1+\gamma)\norm{\nabla_{\bs{\xi}} F(\bs{\xi})}{}^2,
\end{align}
which in turn leads us to
\begin{equation}
\gamma=\dfrac{\bs{\upsilon}\cdot\nabla_{\bs{\xi}}F(\bs{\xi})}{\norm{\nabla_{\bs{\xi}} F(\bs{\xi})}{}^2}-1.
\end{equation}
Then, we may write the explicit form of the error of the gradient, that is,
\begin{equation}
\bs{e}_{\scriptscriptstyle{\perp}}=\bs{\upsilon} - (\bs{\upsilon}\cdot\bs{e}_{\scriptscriptstyle{\nabla}})\bs{e}_{\scriptscriptstyle{\nabla}}\qquad\text{and}\qquad\bs{e}_{\scriptscriptstyle{\parallel}}\coloneqq\gamma\nabla_{\bs{\xi}}F(\bs{\xi}).
\end{equation}

By the Pythagorean theorem, we have that
\begin{align}\label{eq:pointwise.error}
\norm{\bs{\upsilon}-\nabla_{\bs{\xi}}F(\bs{\xi})}{}^2&= \norm{\bs{e}_{\scriptscriptstyle{\parallel}}}{}^2+\norm{\bs{e}_{\scriptscriptstyle{\perp}}}{}^2,\nonumber\\[4pt]
&=\gamma^2\norm{\nabla_{\bs{\xi}}F(\bs{\xi})}{}^2+\norm{\bs{\upsilon} - (\bs{\upsilon}\cdot\bs{e}_{\scriptscriptstyle{\nabla}})\bs{e}_{\scriptscriptstyle{\nabla}}}{}^2\;.
\end{align}
Then, taking the expectation conditioned on $\bs{\xi}$, we arrive at
\begin{align}\label{eq:partwise.error}
\bb{E}[\norm{\bs{\upsilon}-\nabla_{\bs{\xi}}F(\bs{\xi})}{}^2|\bs{\xi}]=\bb{E}[\gamma^2|\bs{\xi}]\norm{\nabla_{\bs{\xi}}F(\bs{\xi})}{}^2+\bb{E}[\norm{\bs{\upsilon} - (\bs{\upsilon}\cdot\bs{e}_{\scriptscriptstyle{\nabla}})\bs{e}_{\scriptscriptstyle{\nabla}}}{}^2|\bs{\xi}].
\end{align}
Note that by Assumption~\ref{as:grad_ass_unbiased}, $\bb{E}[\gamma|\bs{\xi}]=0$. Hence,
\begin{align}\label{eq:gamma.var}
\bb{E}[\gamma^2|\bs{\xi}]&=\bb{V}[\gamma|\bs{\xi}],\nonumber\\[4pt]
&=\bb{V}\bigg[\dfrac{\bs{\upsilon}\cdot\nabla_{\bs{\xi}}F(\bs{\xi})}{\norm{\nabla_{\bs{\xi}} F(\bs{\xi})}{}^2}-1\bigg\vert\bs{\xi}\bigg],\nonumber\\[4pt]
&=\dfrac{1}{\norm{\nabla_{\bs{\xi}} F(\bs{\xi})}{}^4}\bb{V}[\bs{\upsilon}\cdot\nabla_{\bs{\xi}}F(\bs{\xi})|\bs{\xi}].
\end{align}
In view of \eqref{eq:partwise.error} and \eqref{eq:gamma.var}, the statistical error can be written as
\begin{align}\label{eq:statistical.error.decomposition}
\bb{E}[\norm{\bs{\upsilon}-\nabla_{\bs{\xi}}F(\bs{\xi})}{}^2|\bs{\xi}]&=\dfrac{1}{\norm{\nabla_{\bs{\xi}} F(\bs{\xi})}{}^2}\bb{V}[\bs{\upsilon}\cdot\nabla_{\bs{\xi}}F(\bs{\xi})|\bs{\xi}]+\bb{E}[\norm{\bs{\upsilon} - (\bs{\upsilon}\cdot\bs{e}_{\scriptscriptstyle{\nabla}})\bs{e}_{\scriptscriptstyle{\nabla}}}{}^2|\bs{\xi}]\;,\nonumber\\[4pt]
&=\bb{V}[\bs{\upsilon}\cdot\bs{e}_{\scriptscriptstyle{\nabla}}|\bs{\xi}]+\bb{E}[\norm{\bs{\upsilon} - (\bs{\upsilon}\cdot\bs{e}_{\scriptscriptstyle{\nabla}})\bs{e}_{\scriptscriptstyle{\nabla}}}{}^2|\bs{\xi}].
\end{align}
\end{proof}

By Lemma \ref{lm:error.decomposition}, the norm test condition
$\dfrac{\bb{E}[\norm{\bs{\upsilon}_k-\nabla_{\bs{\xi}}
    F(\bs{\xi}_k)}{}^2|\bs{\xi}_k]}{\norm{\nabla_{\bs{\xi}}
    F(\bs{\xi}_k)}{}^2} \le\epsilon^2$ is met if we bound the
statistical estimator such that
\begin{equation}
\left\{
  \begin{aligned}
    \bb{E}[\|\bs{e}_{\scriptscriptstyle{\parallel}}\|^2|\bs{\xi}]&\le\theta^2\|\nabla_{\bs{\xi}}F(\bs{\xi})\|^2,\\[4pt]
    \bb{E}[\|\bs{e}_{\scriptscriptstyle{\perp}}\|^2|\bs{\xi}]&\le\nu^2\|\nabla_{\bs{\xi}}F(\bs{\xi})\|^2,
  \end{aligned}
\right.
\end{equation}
with $\theta^2 + \nu^2 \le \epsilon^2$. The preceding display can be equivalently written as
\begin{equation}\label{eq:conditions}
\left\{
  \begin{aligned}
    \bb{V}[\bs{\upsilon}\cdot\bs{e}_{\scriptscriptstyle{\nabla}}|\bs{\xi}]&\le\theta^2\|\nabla_{\bs{\xi}}F(\bs{\xi})\|^2,\qquad\text{inner product test},\\[4pt]
    \bb{E}[\norm{\bs{\upsilon} - (\bs{\upsilon}\cdot\bs{e}_{\scriptscriptstyle{\nabla}})\bs{e}_{\scriptscriptstyle{\nabla}}}{}^2|\bs{\xi}]&\le\nu^2\|\nabla_{\bs{\xi}}F(\bs{\xi})\|^2,\qquad\text{orthogonality test}.
  \end{aligned}
\right.
\end{equation}
which are the conditions presented by Bollapragada et
al. \cite{Bol18}, for the inner product and orthogonality tests.

To obtain a uniform bound on the variance of the estimator
$\bs{\upsilon}$ at every iteration $k$, we use
$\theta^2+\nu^2=\epsilon^2$. Thus, in view of Proposition~\ref{prp:convergence.mice} the \texttt{SGD} method with
\begin{equation}\label{eq:step.size.strongly.convex.modified}
\eta = \dfrac{2}{(L+\mu)(1 + \theta^2+\nu^2)},
\end{equation}
enjoys linear convergence rate:
\begin{align}\label{eq:residual.recursion.nu.vartheta}
\bb{E}[\norm{\bs{\xi}_{k+1} - \bs{\xi}^*}{}^2] &\le \left[\frac{\left(\dfrac{\kappa-1}{\kappa+1}\right)^2+\theta^2+\nu^2}{1+\theta^2+\nu^2}\right]^k\bb{E}[\norm{\bs{\xi}_0 - \bs{\xi}^*}{}^2]\;.
\end{align}

The convergence rates given in
\eqref{eq:residual.recursion.nu.vartheta} and
\eqref{eq:residual.recursion} are identical for any choice of $\theta$, $\nu$, and
$\epsilon$ if $\theta^2+\nu^2=\epsilon^2$. Thus,
we conclude that controlling the statistical error with
\eqref{eq:conditions} or with \eqref{eq:grad_rel_err.norm}, arising
from identity \eqref{eq:stat.decomposition} in Lemma
\ref{lm:error.decomposition}, are equivalent strategies based on the convergence rates associated with the \texttt{SGD} method.

%% file: section-2.tex
\section{Covariance decomposition}
\label{cov.decomposition}

In this section, we decompose the covariance of the mean gradient
estimator into two components. These components are the covariance in
the orthogonal and in the parallel directions of the true
gradient $\nabla_{\bs{\xi}}F$, respectively. From this covariance
decomposition, we derive formulae to compute the sample sizes to
estimate the mean gradient using two different strategies.

Recall that
\begin{equation}
\norm{\nabla_{\bs{\xi}}F(\bs{\xi})}{}=\nabla_{\bs{\xi}}F(\bs{\xi})\cdot\bs{e}_{\scriptscriptstyle{\nabla}}\;,
\end{equation}
with $\bs{e}_{\scriptscriptstyle{\nabla}}$ as in
Lemma~\ref{lm:error.decomposition}. It then follows from
Assumption~\ref{as:grad_ass_unbiased} that
\begin{align}\label{eq:covariance.grad.direction}
\bb{V}[\bs{\upsilon}\cdot\bs{e}_{\scriptscriptstyle{\nabla}}|\bs{\xi}]&=\bb{E}[\big((\bs{\upsilon}-\nabla_{\bs{\xi}}F(\bs{\xi}))\cdot\bs{e}_{\scriptscriptstyle{\nabla}}\big)^2|\bs{\xi}]\nonumber\\[4pt]
&=\bb{E}[(\bs{\upsilon}-\nabla_{\bs{\xi}}F(\bs{\xi}))\otimes(\bs{\upsilon}-\nabla_{\bs{\xi}}F(\bs{\xi}))|\bs{\xi}]\colon\bs{e}_{\scriptscriptstyle{\nabla}}\otimes\bs{e}_{\scriptscriptstyle{\nabla}}\nonumber\\[4pt]
&=\bs{\Sigma}\colon\bs{e}_{\scriptscriptstyle{\nabla}}\otimes\bs{e}_{\scriptscriptstyle{\nabla}},
\end{align}
where $\bs{\Sigma}\coloneqq\bb{E}[(\bs{\upsilon}-\nabla_{\bs{\xi}}F(\bs{\xi}))\otimes(\bs{\upsilon}-\nabla_{\bs{\xi}}F(\bs{\xi}))|\bs{\xi}]$. Also, note that $\tr\bs{\Sigma}\coloneqq\bs{\Sigma}\colon\id=\bb{E}[\norm{\bs{\upsilon}-\nabla_{\bs{\xi}}F(\bs{\xi})}{}^2|\bs{\xi}]$, where $\id$ is the identity tensor and the symbol $\colon$ represents the double contraction of tensors, that is, $\bs{A}\colon\bs{B}\coloneqq A_{ij}B_ij$ in Einstein's notation.

Next, in view of \eqref{eq:statistical.error.decomposition} along with \eqref{eq:covariance.grad.direction}, we have that
\begin{align}
\bb{E}[\norm{\bs{\upsilon}-\nabla_{\bs{\xi}}F(\bs{\xi})}{}^2|\bs{\xi}]&=\bb{V}[\bs{\upsilon}\cdot\bs{e}_{\scriptscriptstyle{\nabla}}|\bs{\xi}]+\bb{E}[\norm{\bs{\upsilon} - (\bs{\upsilon}\cdot\bs{e}_{\scriptscriptstyle{\nabla}})\bs{e}_{\scriptscriptstyle{\nabla}}}{}^2|\bs{\xi}]\nonumber\\[4pt]
&=\bs{\Sigma}\colon\bs{P}_{\!\scriptscriptstyle{\nabla}}+\bs{\Sigma}\colon\bs{P}_{\!\scriptscriptstyle{\perp}}=\tr\bs{\Sigma},
\end{align}
where $\bs{P}_{\!\scriptscriptstyle{\nabla}}\coloneqq\bs{e}_{\scriptscriptstyle{\nabla}}\otimes\bs{e}_{\scriptscriptstyle{\nabla}}$ and $\bs{P}_{\!\scriptscriptstyle{\perp}}\coloneqq\id-\bs{e}_{\scriptscriptstyle{\nabla}}\otimes\bs{e}_{\scriptscriptstyle{\nabla}}$.

Using this covariance decomposition, the inner product and orthogonality conditions read
\begin{equation}\label{eq:conditions.covariance}
\left\{
  \begin{aligned}
    \bs{\Sigma}\colon\bs{P}_{\!\scriptscriptstyle{\nabla}}&\le\theta^2\|\nabla_{\bs{\xi}}F(\bs{\xi})\|^2,\qquad\text{inner product test}\\[4pt]
    \bs{\Sigma}\colon\bs{P}_{\!\scriptscriptstyle{\perp}}&\le\nu^2\|\nabla_{\bs{\xi}}F(\bs{\xi})\|^2,\qquad\text{orthogonality test}
  \end{aligned}
\right.
\end{equation}

Using the usual Monte Carlo type of statistical estimator
\begin{equation}
\bs{\upsilon}^{(b)}=\dfrac{1}{b}\sum_{\alpha=1}^b\nabla_{\bs{\xi}}f(\bs{\xi},\bs{\vartheta}_\alpha)\;,
\end{equation}
with $\{\bs{\vartheta}_\alpha\}_{\alpha=1,\dots b}$ being an
i.i.d.~sample, we have that the covariance of this statistical
estimator is given by
$\bb{C}ov[\bs{\upsilon}^{(b)}|\bs{\xi}]=\bs{\Sigma}/b$. Moreover, in
that case the conditions in~\eqref{eq:conditions.covariance} are
specialized to
\begin{equation}\label{eq:conditions.inner.orth}
\left\{
  \begin{aligned}
    b_{\scriptscriptstyle{\nabla}}&\ge\dfrac{\bs{\Sigma}\colon\bs{P}_{\!\scriptscriptstyle{\nabla}}}{\theta^2\|\nabla_{\bs{\xi}}F(\bs{\xi})\|^2},\qquad\text{inner product test}\\[4pt]
    b_{\scriptscriptstyle{\perp}}&\ge\dfrac{\bs{\Sigma}\colon\bs{P}_{\!\scriptscriptstyle{\perp}}}{\nu^2\|\nabla_{\bs{\xi}}F(\bs{\xi})\|^2},\qquad\text{orthogonality test}
  \end{aligned}
\right.
\end{equation}
We refer to condition \eqref{eq:conditions.inner.orth} with
$b\coloneqq\max(b_{\scriptscriptstyle{\nabla}},b_{\scriptscriptstyle{\perp}})$
as inner/orth test sample size.

Conversely, the norm test \eqref{eq:grad_rel_err.norm}, given in
Definition \ref{df:grad_rel_err.norm}, implies the sample size
condition
\begin{equation}\label{eq:conditions.norm}
\bar{b}\ge\dfrac{\bs{\Sigma}\colon\id}{\epsilon^2\|\nabla_{\bs{\xi}}F(\bs{\xi})\|^2}.
\end{equation}
We referred to $\bar{b}$ in condition \eqref{eq:conditions.norm} as
the norm test sample size.

%% file: section-3.tex
\section{Complexity}
\label{complexity}

In this section, we demonstrate that the inner/orth test's algorithmic complexity, in terms of
gradient evaluations at a given iteration, is equal to the norm test's algorithmic complexity provided that
$\epsilon^2=\theta^2+\nu^2$ and $\theta$ and $\nu$ are optimally chosen, otherwise the norm test is always more computational affordable. In what follows, we measure the complexity in terms of
gradient evaluations, which are determined by the sample sizes.

Our point of departure is that
\begin{equation}\label{eq:param.relation}
\epsilon^2=\theta^2+\nu^2.
\end{equation}
Note that $\epsilon=1$ is asymptotically optimal for
strongly-convex functions, as is detailed in~\cite[Remark 10]{Car21}.
In general, given a fixed $\epsilon$ we have that if we increase
$\theta$, $\nu$ should be decreased for a fixed $\epsilon>0$ as per relation \eqref{eq:param.relation}. Bearing this in mind and recalling the inner/orth test sample size
$b\coloneqq\max(b_{\scriptscriptstyle{\nabla}},b_{\scriptscriptstyle{\perp}})$,
the optimal batch size $b^\ast$ is obtained when both constraints
\eqref{eq:conditions.covariance} are active, namely
\begin{equation}\label{eq:optimal.b}
b^\ast=b_{\scriptscriptstyle{\nabla}}=b_{\scriptscriptstyle{\perp}}.
\end{equation}
Thus, for the relation $b_{\scriptscriptstyle{\nabla}}=b_{\scriptscriptstyle{\perp}}$ to be true, we have that
\begin{align}\label{eq:identity}
\dfrac{\theta^2}{\nu^2}&=\dfrac{\bs{\Sigma}\colon\bs{P}_{\!\scriptscriptstyle{\nabla}}}{\bs{\Sigma}\colon\bs{P}_{\!\scriptscriptstyle{\perp}}}\nonumber\\[4pt]
&=\dfrac{\bs{\Sigma}\colon\id-\bs{\Sigma}\colon\bs{P}_{\!\scriptscriptstyle{\perp}}}{\bs{\Sigma}\colon\bs{P}_{\!\scriptscriptstyle{\perp}}}\nonumber\\[4pt]
&=\dfrac{\bs{\Sigma}\colon\id}{\bs{\Sigma}\colon\bs{P}_{\!\scriptscriptstyle{\perp}}}-1.
\end{align}
Moreover, using expression \eqref{eq:param.relation} along with \eqref{eq:identity} we arrive at
\begin{equation}
1+\dfrac{\theta^2}{\nu^2}=\dfrac{\epsilon^2}{\nu^2}=\dfrac{\bs{\Sigma}\colon\id}{\bs{\Sigma}\colon\bs{P}_{\!\scriptscriptstyle{\perp}}}.
\end{equation}
From which, we are led to
\begin{equation}\label{eq:nu.epsilon.relation}
\nu^2=\epsilon^2\dfrac{\bs{\Sigma}\colon\bs{P}_{\!\scriptscriptstyle{\perp}}}{\bs{\Sigma}\colon\id}.
\end{equation}
Analogously, we have that
\begin{equation}\label{eq:vartheta.epsilon.relation}
\theta^2=\epsilon^2\dfrac{\bs{\Sigma}\colon\bs{P}_{\!\scriptscriptstyle{\nabla}}}{\bs{\Sigma}\colon\id}.
\end{equation}
Then, in view of \eqref{eq:nu.epsilon.relation} along with \eqref{eq:optimal.b}, we have that
\begin{equation}\label{eq:optimal.b.explicit}
b^\ast=b_{\scriptscriptstyle{\nabla}}=b_{\scriptscriptstyle{\perp}}=\dfrac{\bs{\Sigma}\colon\id}{\epsilon^2\|\nabla_{\bs{\xi}}F(\bs{\xi})\|^2}.
\end{equation}
Notice that the optimal inner/orth test samples size $b^\ast$ identified
above is identical to the norm test sample $\bar{b}$
in~\eqref{eq:conditions.norm}.

Now, we demonstrated that both strategies, namely the norm and the inner
product/orthogonality tests, are also equivalent in terms of
computational cost if $\nu$ is chosen according to
\eqref{eq:nu.epsilon.relation} and $\theta$ according to
\eqref{eq:vartheta.epsilon.relation}. Finally, we emphasize that this
equivalence holds in the asymptotic regime, given that the convergence
results assume $k\to\infty$.

%% file: section-4.tex
\section{Numerical experiments}
\label{numerics}

We compare the inner product/orthogonality tests and the norm test for
two different the objective functions considering $\theta$, $\nu$, and $\epsilon$ fixed throughout the optimization while the relation $\epsilon^2=\theta^2+\nu^2$ holds.

\noindent \emph{Objective function 1 in $\bb{R}^3$:}
\begin{equation}\label{eq:obj.F.1}
\bb{E} [f(\bs{\xi},\bs{\vartheta})|\bs{\xi}] = \bb{E}\left[\dfrac{1}{2}\bs{\xi}\cdot\bs{H}\bs{\xi}-\bs{\vartheta}\cdot\bs{\xi} | \bs{\xi}\right],\qquad
\bs{H}=
\begin{bmatrix}
2 & 1 & 1 \\
1 & 10 & 1 \\
1 & 1 & 100 \\
\end{bmatrix},
\qquad\bs{\vartheta}\sim\cl{N}(\bs{0},10^3\id_{3\times{3}}),
\end{equation}
where $\id_{3\times{3}}$ is the identity of dimension $3\times{3}$. The optimal point of this problem is $\bs{\xi}^* = (0,0,0)$.

\noindent \emph{Objective function 2 in $\bb{R}^2$:}
\begin{equation}
\bb{E}[f(\bs{\xi}, \vartheta)|\bs{\xi}] = \bb{E}[\frac{1}{2} \bs{\xi} \cdot \bs{H}(\vartheta) \bs{\xi}|\bs{\xi}] - \bs{b} \cdot \bs{\xi},\qquad
\bs{H}(\vartheta) \coloneqq
\id_{2\times{2}}(1-\vartheta) +
\begin{bmatrix}
 2 \kappa & 0.5 \\
 0.5      & 1
\end{bmatrix}
\vartheta,
\end{equation}
where $\id_{2\times{2}}$ is the identity of dimension $2\times{2}$,
$\bs{b}$ is a vector of ones, and $\vartheta \sim \cl{U}(0, 1)$. Here,
we use $\kappa=100$.  That is, the objective function to be minimized
is
\begin{equation}\label{eq:obj.F.2}
  F(\bs{\xi}) = \frac{1}{2} \bs{\xi} \cdot
  \bb{E}[\bs{H}(\bs{\vartheta})] \bs{\xi}
  - \bs{b} \cdot \bs{\xi},
\end{equation}
where
\begin{equation}
 \bb{E}[\bs{H}(\vartheta)] =
 \begin{bmatrix}
  \kappa+0.5 & 0.25 \\
  0.25       & 1
 \end{bmatrix}.
\end{equation}
The optimal point of this problem is
$\bs{\xi}^* = \bb{E}[\bs{H}(\theta)]^{-1} \bs{b}$.

\begin{figure}
  \subfloat[Cases: $\#1$, $\#2$, and $\#3$]{\includegraphics[width=0.45\textwidth]{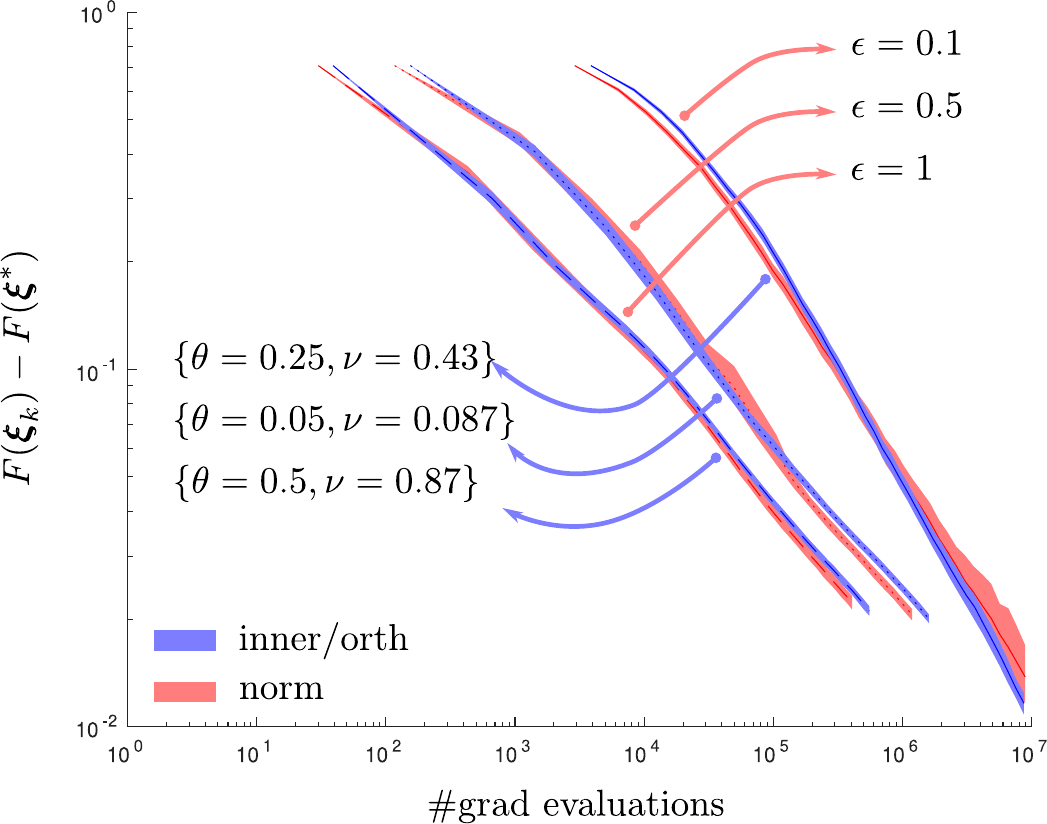}\label{fg:1.a}}
  \subfloat[Cases: $\#3$ and $\#4$]{\includegraphics[width=0.45\textwidth]{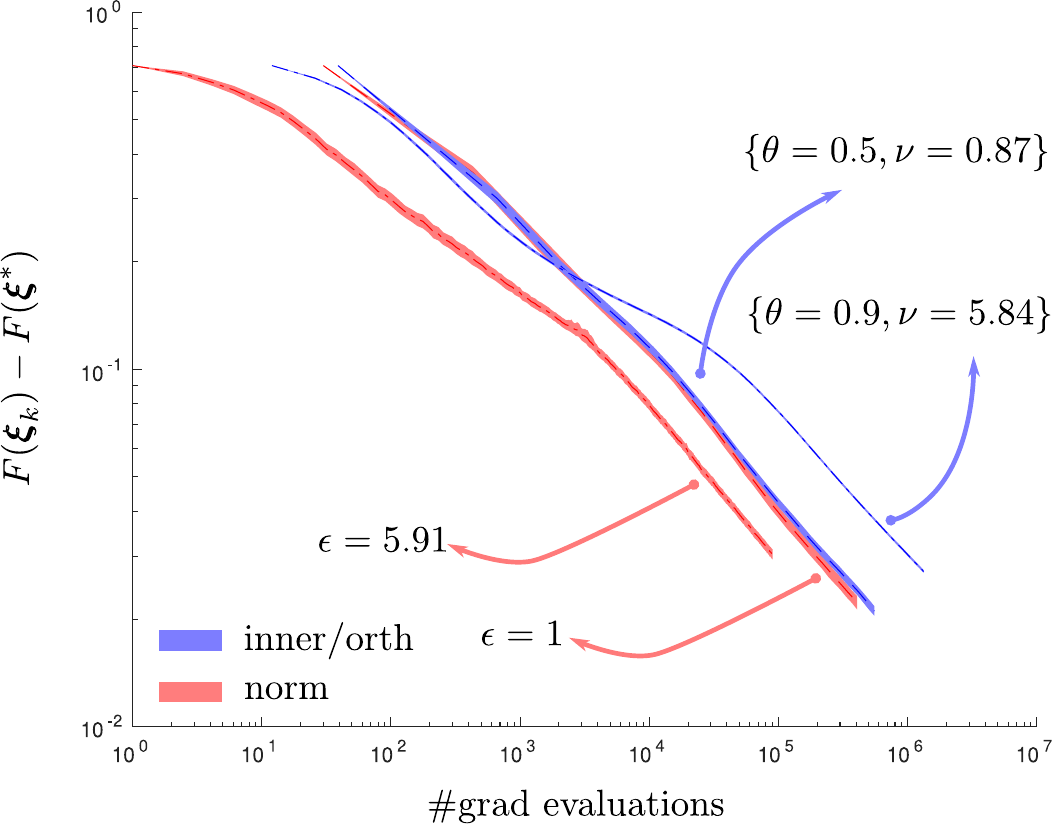}\label{fg:1.b}}
  \caption{The objective function in \eqref{eq:obj.F.1}: optimality gap versus total cost.}
  \label{fg:gapxcost.quad.1}
\end{figure}

\begin{figure}
  \subfloat[Cases: $\#1$, $\#2$, and $\#3$]{\includegraphics[width=0.45\textwidth]{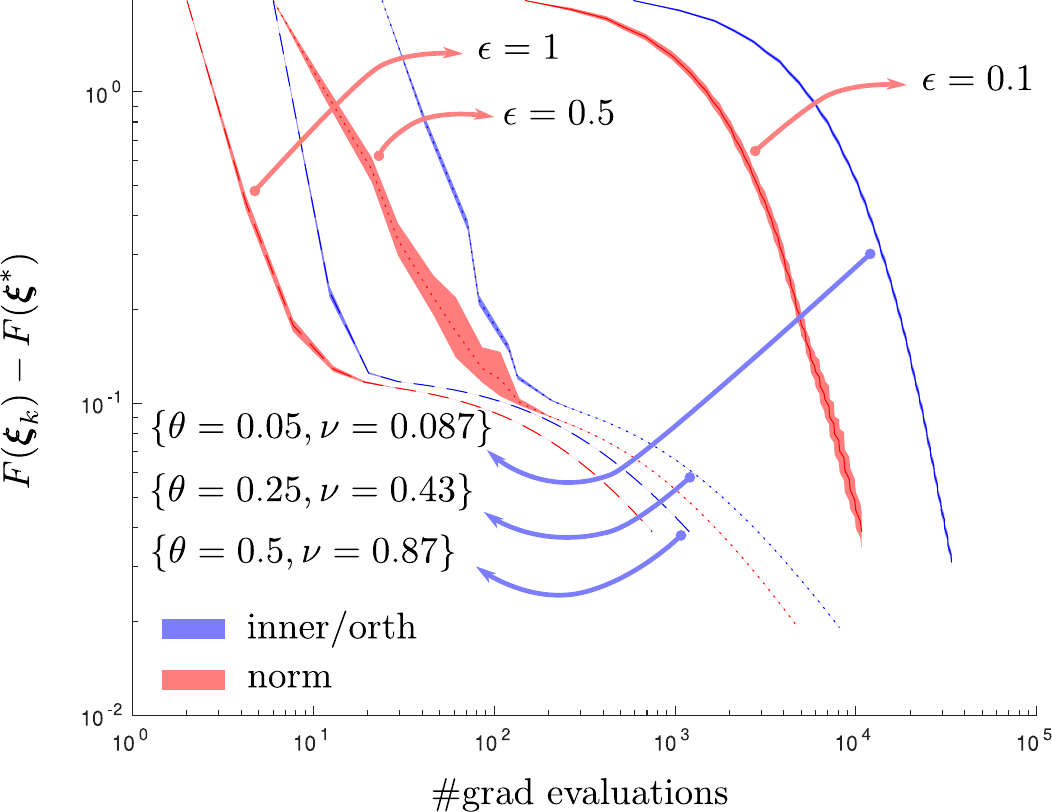}\label{fg:2.a}}
  \subfloat[Cases: $\#3$ and $\#4$]{\includegraphics[width=0.45\textwidth]{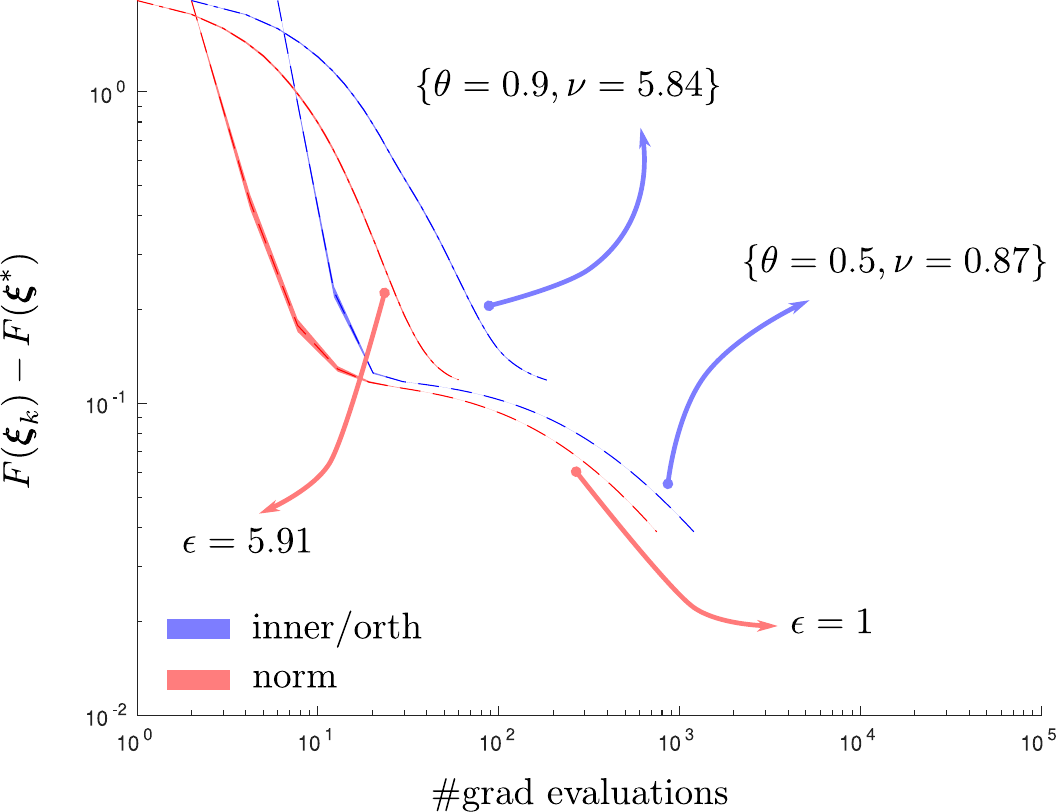}\label{fg:2.b}}
  \caption{The objective function in \eqref{eq:obj.F.2}: optimality gap versus total cost.}
  \label{fg:gapxcost.quad.2}
\end{figure}

Table \ref{tb:parameters} shows the set of parameters
used to compare four different scenarios. Since we keep $\theta$ and $\nu$ constant throughout the optimization, we cannot expect a complexity equivalence, however, the convergence rates are still equivalent. In terms of convergence, the first three scenarios show the effect of increasing the relative tolerance in the statistical estimator of the mean gradient. In the last scenario, $\#4$, we use the parameters $\theta=0.9$ and $\nu=5.84$ found in \cite{Bol18}, and compare them against $\epsilon=5.91$. In Figures
\ref{fg:gapxcost.quad.1} and \ref{fg:gapxcost.quad.2}, we present the
optimality gap $F(\bs{\xi}_k)-F(\bs{\xi}^\ast)$ versus the total cost
(number of gradient evaluations), namely $\sum_{k=0}^K b_k$ for the
cases in Table \ref{tb:parameters}. To generate this plot, we solve
the optimization problems $1000$ times, each repetition is initialized
with the same $\bs{\xi}_0=[0.225, -0.2, 0.1]$ for objective function 1
($\bs{\xi}_0=[20.0, 50.0]$ for objective function 2) but has been
performed with different random seeds. We then present the $95\%$ confidence
interval of these runs. In these convergence
plots, we show the equivalence in terms of cost and
convergence rates as theoretically established in previous sections
regardless of the choices of $\epsilon$ for the norm test, and $\theta$
and $\nu$ for the inner/orth test as long as the relation
$\epsilon^2=\theta^2+\nu^2$ holds.
\begin{table}
\caption{Parameters used in the simulations, satisfying $\epsilon^2=\theta^2+\nu^2$.}
\begin{tabular}{ l | l | l | l }
case & $\epsilon$ & $\theta$ & $\nu$ \\\hline
$\#1$& 0.1 & 0.05 & 0.087 \\
$\#2$& 0.5 & 0.25 & 0.43 \\
$\#3$& 1 & 0.5 & 0.87 \\
$\#4$& 5.91 & 0.9 & 5.84 \\
\end{tabular}
\label{tb:parameters}
\end{table}
The equivalence in terms of both convergence and cost is clear from
Figure \ref{fg:gapxcost.quad.1}\subref{fg:1.a} for the first objective
function \eqref{eq:obj.F.1} for the parameter cases $\#1$, $\#2$, and
$\#3$.  The equivalence does not hold, however, for the case
$\#4$. Furthermore, from Figure \ref{fg:gapxcost.quad.2} we deduce
that the equivalence does not hold in any of the considered cases for
objective function \eqref{eq:obj.F.2}. This is because, for the
equivalence between these batch size selection approaches, both
expressions \eqref{eq:nu.epsilon.relation} and
\eqref{eq:vartheta.epsilon.relation} must hold at every iteration of
the \texttt{SGD} method. To achieve this, and hence be equivalent
selection strategies, the parameters $\theta$ and $\nu$ should thus
vary along the optimization path.

\section{Acknowledgments}

This work was partially supported by the KAUST Office of Sponsored Research (OSR) under Award numbers URF$/1/2281-01-01$, URF$/1/2584-01-01$ in the KAUST Competitive Research Grants Program Round 8, the Alexander von Humboldt Foundation.

\section{Conclusions}

In this study, we demonstrate that the norm and inner product/orthogonality tests are theoretically equivalent in terms of cost and convergence rates. We demonstrate the equivalence of our theoretical predictions by illustrating them in a simple stochastic optimization problem.

\section{Data availability statement}

The data that support the findings of this study are available from the corresponding author upon resonable request.

\subsection{Conflict of interest}

The authors have no conflicts to disclose.

%% file: main.bbl
\begin{thebibliography}{1}

\bibitem{Bol18}
Raghu Bollapragada, Richard Byrd, and Jorge Nocedal.
\newblock Adaptive sampling strategies for stochastic optimization.
\newblock {\em SIAM Journal on Optimization}, 28(4):3312--3343, 2018.

\bibitem{Car91}
Richard~G Carter.
\newblock On the global convergence of trust region algorithms using inexact
  gradient information.
\newblock {\em SIAM Journal on Numerical Analysis}, 28(1):251--265, 1991.

\bibitem{Byr12}
Richard~H Byrd, Gillian~M Chin, Jorge Nocedal, and Yuchen Wu.
\newblock Sample size selection in optimization methods for machine learning.
\newblock {\em Mathematical programming}, 134(1):127--155, 2012.

\bibitem{Has14}
Fatemeh~S Hashemi, Soumyadip Ghosh, and Raghu Pasupathy.
\newblock On adaptive sampling rules for stochastic recursions.
\newblock In {\em Proceedings of the Winter Simulation Conference 2014}, pages
  3959--3970. IEEE, 2014.

\bibitem{Car18}
Coralia Cartis and Katya Scheinberg.
\newblock Global convergence rate analysis of unconstrained optimization
  methods based on probabilistic models.
\newblock {\em Mathematical Programming}, 169(2):337--375, 2018.

\bibitem{Car21}
Andre Carlon, Luis Espath, Rafael Lopez, and Raul Tempone.
\newblock Multi-iteration stochastic optimizers.
\newblock {\em arXiv preprint arXiv:2011.01718}, 2020.

\end{thebibliography}
